\documentclass[10pt]{article}
\usepackage{amsmath,amssymb,amsthm,amscd}
\numberwithin{equation}{section}

 \DeclareMathOperator{\Aut}{Aut}
 \DeclareMathOperator{\Reg}{Reg}

\def\ddc{\sqrt{-1}\partial{\bar\partial}}

\def\CC{{\mathbb C}}

\newtheorem{prop}{Proposition}[section]
\newtheorem{theo}[prop]{Theorem}
\newtheorem{lemm}[prop]{Lemma}
\newtheorem{coro}[prop]{Corollary}
\newtheorem{rema}[prop]{Remark}
\newtheorem{exam}[prop]{Example}
\newtheorem{defi}[prop]{Definition}

\def\begeq{\begin{equation}}
\def\endeq{\end{equation}}

\def\and{\quad{\rm and}\quad}

\def\<{\langle}
\def\>{\rangle}

\def\om{\omega}
\def\al{\alpha}
\def\Om{\Omega}

\def\lbr{\lbrace}
\def\rbr{\rbrace}

\begin{document}

\title{Stability of Bounded Solutions for Degenerate 
Complex Monge-Amp\`ere Equations}

\author{S\l awomir Dinew\\
Jagiellonian University, Krak\'ow\\
Zhou Zhang\\ 
Department of Mathematics,\\
University of Michigan, at Ann Arbor}
\date{}
\maketitle

\begin{abstract}
We show a stability estimate for the degenerate complex 
Monge-Amp\`ere operator that generalizes a result of 
Ko\l odziej \cite{koj}. In particular, we obtain the optimal 
stability exponent and also treat the case when the right 
hand side is a general Borel measure satisfying certain 
regularity conditions. Moreover our result holds for functions 
plurisubharmonic with respect to a big form generalizing 
thus the K\"ahler form setting in \cite{koj}.
\end{abstract}

\section{Introduction and the Main Theorem}

In this work, we generalize and strengthen Ko\l odziej's 
stability result concerning bounded solutions for 
complex Monge-Amp\`ere equations, which is 
summarized in \cite{koj} (see also \cite{kojnotes}). 
The solutions are understood in the sense of 
pluripotential theory, i.e. we do not impose any other
regularity than upper semicontinuity and 
boundedness. It is, however, a classical fact that the
image of the Monge-Amp\`ere operator can be well 
defined as a Borel measure in this setting.    

The equation is considered over a closed K\"ahler 
manifold $X$\ of complex dimension $n\geqslant 2
$ \footnote{When $n=1$, the manifold is a Riemann 
surface and Monge-Amp\`ere operator is just the 
Laplace operator.}.   

Suppose $\omega$ is a real smooth closed semi-positive 
$(1,1)$-form over $X$, $\Omega$ is a positive Borel 
measure on $X$ and $f\in L^p(X)$ for some $p>1$ 
is non-negative, where the definition of the function space 
$L^p(X)$ is with respect to $\Omega$. The equation we 
consider is  
$$(\omega+\sqrt{-1}\partial\bar{\partial }u)^n=f\Omega.$$  
Using $d=\partial+\bar{\partial }$\ and $d^c:=\frac{\sqrt{-1}}{2}
(\bar{\partial}-\partial)$ we have $dd^c=\sqrt{-1}\partial\bar
{\partial }$ and this convention is also often used in the 
literature.

As mentioned above, we require regularity of $u$ 
much less than what is needed to make pointwise 
sense for the left hand side. More specifically, we 
look for solutions in the function class $PSH_\omega(X)\cap 
L^\infty(X)$, where $u\in PSH_\omega(X)$ means that $\omega+\sqrt{-1}\partial
\bar{\partial }u$ is non-negative in the sense of 
distribution theory. 

Of course, there is an obvious condition for the 
existence of such a solution coming from global 
integration over $X$, i.e. $\int_X{\omega }^n=\int_Xf\Omega$. 
This condition follows from Stokes theorem in the smooth case, and hence (by smooth 
approximation) in our case either.

Ko\l odziej mainly studied the case when $\omega$ is a 
K\"ahler metric, or equivalently, $[\omega]$ is a K\"ahler 
class, and $\Omega$\ is a smooth volume form. The 
existence of bounded solution in this 
case is proved. In fact, even more general 
$f$'s than $L^p$ functions are treated in \cite{koj98}, but for our main 
concern, we restrict to $L^p$ functions. Further, 
in this case, the bounded solution is always 
continuous as proved in \cite{koj98}. So 
in the discussion of stability there, continuity of the 
solutions is naturally assumed.  

The degeneration we want to consider in this note is 
in two places.

First we allow $\omega$\ to be just semi-positive instead of 
being K\"ahler; we are especially interested in the case 
when $\omega$ is the pullback of a K\"ahler metric under 
a holomorphic map preserving dimensions. The 
following theorem from \cite{zh} gives the precise 
picture of $\omega$ and the corresponding existence result. 
This result uses an argument very close to Ko\l odziej's. 
Both of them have found the notion of relative capacity, 
introduced in \cite{bed-tay}, extremely useful. 

\begin{theo}\label{1.1}

Let $X$ be a closed K\"ahler manifold with 
(complex) dimension $n\geqslant 2$. 
Suppose we have a holomorphic map $F: X
\to\mathbb{CP}^N$ with the image $F(X)$ of 
the same dimension as $X$. Let $\omega_M$ be 
any K\"ahler form over some neighbourhood 
of $F(X)$ in $\mathbb{CP}^N$. For the 
following equation of Monge-Amp\`ere type:  
\begin{equation}
(\omega+\sqrt{-1}\partial\bar\partial u)^n=f\Omega,\nonumber 
\end{equation}
where $\omega=F^*\omega_M$, $\Omega$ is a fixed smooth 
(non-degenerate) volume form over $X$ and 
$f$ is a nonnegative function in $L^p(X)$ for 
some $p>1$ with the correct total integral 
over $X$, i.e. $\int_X f\Omega=\int_X(F^*\omega_M)^n$, then 
we have the following:  

(1) (A priori estimate) If $u$ is a weak solution 
in $PSH_{\omega}(X)\cap L^{\infty}(X)$ of the equation with 
the normalization $sup_X u=0$, then there is 
a constant $C$ such that $\|u\|_{L^\infty}\le C\|f\|^n_{
L^p}$ where $C$ only depends on $F$, $\omega$ 
and $p$; 

(2) (Existence of a bounded solution) There 
exists a bounded (weak) solution for this 
equation;  

(3) (Continuity and uniqueness of bounded 
solution) If $F$ is locally birational, any 
bounded solution is actually the unique 
continuous solution.

\end{theo}

The a priori estimate was obtained independently in \cite{EGZ} 
(even for more general big forms), and later generalized to more 
singular right hand side in \cite{DP}. As for the continuity of the 
solution, despite serious effort, the situation is still a little bit unclear. 
It is not known whether continuity holds when $\omega$\ is a 
general semi-positive closed form with continuous (even smooth) 
potentials and positive total integral. This problem has attracted 
much interest recently, and for this reason we take the opportunity 
to present a detailed proof of the continuity in the situation above 
Indeed, the argument in \cite{zh} \footnote{See also in \cite{thesis} 
where it is not so easily separated from the context.} is a bit too 
sketchy therefore hard to follow. This will be done in Section 
\ref{continuity}.

 Regardless of that in our discussion of stability we do not impose 
 {\it a priori} continuity of the solutions. The methods we use are 
 independent of that assumption. So, theoretically, solutions might 
 be discontinuous in general, but uniformly close to each other if 
 we perturb the data a little. Needless to say, this is quite an artificial 
 situation. So our results strongly support (but in no way prove) 
 the common belief that continuity holds in general.    

Our second degeneration is that we allow 
$\Omega$\ on the right hand side to be a Borel 
measure instead of smooth volume form. 
Then some restrictions must be imposed, 
since weak solutions for such an equation 
might not be bounded anymore (for example, 
if $\Omega$\ is the Dirac delta measure at some 
point). Worse yet, there are measures for 
which existence of solutions (bounded or 
not) is not known so far. Therefore we 
impose some seemingly natural conditions 
on $\Omega$\ that guarantee boundedness of the solutions. 

\begin{defi} 

We say that a Borel measure is {\it well 
dominated by capacity for $L^p$\ functions}, 
if there exist constants $\alpha>0$\ and $\chi>0$, such that 
for any compact $K\subset X$\ and any non-negative 
$f\in L^p(\Omega),\ p>1$ one has for some constant 
$C$\ independent of $K$, (but dependent on $f$)
$$\Om(K)\leq C cap_{\om}(K)^{1+\alpha},\ \ \and \int_Kf
\Om\leq C cap_{\om}(K)^{1+\chi}$$

\end{defi}

A very similar notion (only the first condition is 
imposed) is discussed in \cite{EGZ}. Both are 
variations of the so-called condition (A), introduced 
by S. Ko\l odziej in \cite{koj98}. These conditions 
(which actually are stronger than condition(A)) 
force boundedness for the solutions $u$\ of
$$(\omega+\sqrt{-1}\partial\bar\partial u)^n=f\Omega$$ 
(see \cite{koj98} for the case $\om$\ is K\"ahler, 
and \cite{EGZ} for the case $\om$\ is merely 
semi-positive).

A few words on the second assumption. When $\Omega$\ is a smooth volume form it is known (again see \cite{koj98} and \cite{EGZ}) that the first condition is satisfied for every $\alpha>0$. Hence by an elementary application of the H\"older inequality the second condition is also satisfied (for every $\chi>0$). The same reasoning also shows that the second condition is a consequence of the first provided $p$\ is big enough (if $\frac{(1+\alpha)(p-1)}{p}>1$). Anyway, one has to impose some condition, since {\it a priori} $f\Om$\ is more singular than $\Om$. 

Note that, as in \cite{koj98} or \cite{kojnotes} the exponent $\chi>0$\ is used to construct an {\it admissible} function $Q$\ with proper polynomial growth and afterwards function $\kappa$\ and its inverse $\gamma$\ (see below for a discussion). When the volume form is smooth one can take arbitrary $\chi>0$\ (of course the bigger $\chi$\ we take, the better). Using this, in \cite{koj} it was shown that one can produce a function $\gamma(t)$\ with growth like $t^{\epsilon},\ \forall \epsilon>0$\ near $0$. When $\chi$\ is bounded from above (i.e. we assume it is a fixed constant dependent on the measure $\mu$), calculations as in \cite{kojnotes} or \cite{koj} show that one can take $\gamma(t)\approx t^{\frac n{\chi}}$. In order to avoid too much technicalities throughout the note we shall work with the assumption that $\chi$\ can be taken arbitrarily large. At the end (see Remark \ref{rem4.1}) we will explain how to modify the argument in the case of fixed $\chi$ and obtain the stability exponent in this case either.

As mentioned in the thesis of the second 
named author \cite{thesis}, Ko\l odziej's 
original argument is almost good enough 
for us except for two issues. One of them, 
about Comparison Principle, is doable 
using the regularizing result in \cite{blo-koj}. 
The other one, an inequality for mixed 
Monge-Amp\`ere measures, looks hard to 
justify for bounded functions. Recently, this 
has been treated by the first named author in 
\cite{dinew} for even more general class of 
functions. 

Now let's state the main theorem. 

\begin{theo}

In the same set-up as in the theorem above (we assume that $\Om$\ is well dominated by capacity for $L^p$\ functions), for 
any non-negative $L^{p}(\Om)$-functions $f$ and $g$ 
with $p>1$ which have the proper total integral 
over $X$, i.e., $\int_X f\Omega=\int_X g\Omega=\int_X {\omega}
^n$, suppose that $\phi$ and $\psi$ in $PSH_{\omega}\cap L
^\infty(X)$ satisfy ${\omega_{\phi}}^{n}=f{\omega}^{n}$ and ${\omega_{\psi}
}^{n}=g {\omega}^{n}$ respectively and are normalized by the conditions $max_{X}\{\phi-\psi\}=max_{X}\{\psi-\phi\}$. Let also $\epsilon>0$\ be arbitrary.

If $\|f-g\|_{L^{1}}\leqslant\gamma (t)t^{n+\epsilon}$ for $\gamma(t)=C\kappa
^{-1}(t)$ with some proper non-negative constant 
$C$ depending only \footnote{The manifold $X$ 
and  metric $\omega$\ also affect $C$.} 
on the $L^{p}$-norms of $f$ and $g$, 
where $\kappa^{-1}(t)$ the inverse function of  
the following $\kappa$ function, 
$$\kappa(r)=C_{n}A^{\frac{1}{n}}\bigl(\int_{r^{-\frac{1}{n}}}^{\infty}
y^{-1}(Q(y))^{-\frac{1}{n}}dy+\bigl(Q(r^{-\frac{1}{n}})
\bigr)^{-\frac{1}{n}}\bigr),$$ 
where $C_{n}$ is a positive constant only depending 
on the complex dimension $n$ and $Q$ is an 
increasing positive function with proper polynomial 
growth, then we can conclude that 
$$\|\phi-\psi\|_{L^{\infty}}\leqslant Ct$$ 
for $t<t_{0}$ where $t_{0}>0$ depends on $\gamma$ and 
$C$ depends on the $L^{p}$-norms of $f$ and $g$.

\end{theo}

As a direct application, we have uniqueness of bounded 
solution from Theorem $(1.1)$. 

Another corollary is the following stability estimate.

\begin{coro}\label{exp2}

In the same setting as above there exists a constant $c=c(p,\epsilon,c_0)$\ where $c_0$\ is an upper bound for $||f||_p$ and $||g||_p$\ such that
$$||\phi-\psi||_{\infty}\leq c||f-g||_1^{\frac 1{n+\epsilon}}$$

\end{coro}

\begin{rema} 

The exponent in the last corollary is improved compared to \cite{koj}. As example \ref{exa} shows, the exponent we obtain is optimal.

\end{rema}
\begin{rema} The Monge-Amp\`ere equation with $\omega$\ big instead of K\"ahler has been studied extensively in the recent years (see \cite{BGZ}, \cite{DP}, \cite{EGZ}). 
\end{rema}

The applications of the result above could go in two directions. The semi-positivity is particularly interesting in geometry, since  the situation we have described above appears naturally in the study of algebraic manifolds of general type (or big line bundles in general) (see e.g. \cite{t-znote}). The degeneration of the measure on the right hand side, in turn, might be useful in complex dynamics and pluripotential theory. Complex dynamics often deals with such singular measures and it is an important question to obtain any regularity for the potential of such measures. The same question is crucial in pluripotential theory while studying extremal functions.  

{\bf Acknowledgment.} 
The authors would like to thank professor S. Ko\l odziej for all the generous help in the formation of this work and beyond. His suggestion for such a joint work is also very important for beginners like us. This work was initiated during the second named author's visit at MRSI (Mathematical Sciences Research Institute) and he would like to thank the institute and the department of Mathematics at University of Michigan, at Ann Arbor, for the arrangement to provide such a wonderful opportunity.

\section{Stability for Nondegenerate Monge-Amp\`ere 
Equations}

For readers' convenience, Ko\l odziej's stability argument 
will be included here. We are going to use global version 
of the notions, for example, capacity for the closed 
manifold $X$. 

Specifically, in this part all the plurisubharmonic 
functions with respect to the K\"ahler metric $\omega$ 
($\omega$-$PSH$ for short) are continuous by definition. 
As explained before, this brings no difference in 
this case. So Comparison Principle between them 
can be justified by the Richberg's approximation 
as in \cite{koj}. 

Basically, all the following argument is directly quoted 
from \cite{koj}.\\ 

Claim: Let $\phi, \psi\in PSH_\omega(X)$\ and satisfy $0\leqslant\phi\leqslant 
C$, then for $s<C+1$, we have 
$$Cap_{\omega}(\{\psi+2s<\phi\})\leqslant\bigl(\frac{C+1}{s}\bigr) 
^{n}\int_{\{\psi+s<\phi\}}(\omega+\sqrt{-1}\partial\bar{\partial}
\psi)^{n}.$$ 

\begin{proof}

Define $E(s):=\{\psi+s<\phi\}$. Take any $\rho\in PSH_{\omega}(X)$ 
valued in $[-1, 0]$. Set $V=\{\psi <\frac{s}{C+1}\rho+(1-
\frac{s}{C+1})\phi-s\}$. 
Since $-s\leqslant\frac{s}{C+1}\rho-\frac{s}{C+1}\phi
\leqslant 0$, we can easily deduce the following 
chain relation of sets: 
$$E(2s)\subset V\subset E(s).$$ 
Then we can have the following computation (with notation  
$\omega_\rho:=\omega+\sqrt{-1}\partial\bar{\partial}\rho$):

\begin{equation}
\begin{split}
(\frac{s}{C+1})^{n}\int_{E(2s)}(\omega+\sqrt{-1}\partial\bar{
\partial}\rho)^{n} 
&\leqslant \int_{V}(\frac{s}{C+1}\omega_{\rho}+(1-\frac{s}{C+
1})\omega_{\phi})^{n}\\
&\leqslant \int_{V}{\omega_{\psi}}^{n} \leqslant \int_{E(s)}{\omega_{\psi}}
^{n}\nonumber
\end{split}
\end{equation}
by the relation of sets above and applying 
comparison principle for the two functions 
appearing in the definition of the set $V$. 

Finally we can conclude the result from the 
definition of $Cap_{\omega}$. 

\end{proof}

Now we state the following version of stability 
result, which is slightly weaker than the result in \cite{kojnotes}. 

\begin{theo}

In the same set-up as before, for any nonnegative 
$L^{p}$-functions $f$ and $g$ with $p>1$ which 
have the proper total integral over $X$, i.e., $\int_X 
f{\omega}^n=\int_X g{\omega}^n=\int_X {\omega}^n$, suppose that $\phi$ 
and $\psi$ in $PSH_{\omega}(X)$ satisfy ${\omega_{\phi}}^{n}=f{\omega}
^{n}$ and ${\omega_{\psi}}^{n}=g {\omega}^{n}$ respectively and 
are normalized by the condition $max_{X}\{\phi-\psi\}=max_{X}\{\psi
-\phi\}$. 

If $\|f-g\|_{L^{1}}\leqslant\gamma (t)t^{n+3}$ for $\gamma(t)=C\kappa
^{-1}(t)$ with some proper nonnegative constant 
$C$ depending only on the $L^{p}$-norms of $f$ 
and $g$ \footnote{The dependence on the manifold 
$X$ and K\"ahler metric $\omega$ should be clear.} , 
where $\kappa^{-1}(t)$ the inverse function of 
the $\kappa$ function in the main theorem, then we can 
conclude that 
$$\|\phi-\psi\|_{L^{\infty}}\leqslant Ct$$ 
for $t<t_{0}$ where $t_{0}>0$ depends on $\gamma$ and 
$C$ depends on the $L^{p}$-norms of $f$ and $g$.

\end{theo}

\begin{proof}

Suppose $\|f\|_{L^{p}}, \|g\|_{L^{p}}\leqslant A$. We will 
be careful about the fact that the constants in the 
argument will only depend on $A$ and the function 
$\gamma$. 

For simplicity, let us normalize to have $\int_{X}{\omega}^{n}=1$. 
And in fact, we can also assume $max_{X}\{\phi-\psi\}=max_{X}\{
\psi-\phi\}>0$ since the case for $=0$ is trivial 
\footnote{In this case, we can have $\phi-\psi\leqslant 0$ and $
\psi-\phi\leqslant 0$, which says $\phi=\psi$. In other words, we 
have the compatible direction.}. 

Without loss of generality, assume $\int_{\{\psi<\phi\}}(f+g)
{\omega}^{n}\leqslant 1$, since $\int_{X}f{\omega}^{n}=\int_{X}
g{\omega}^{n}=1$ and, if needed, one can interchange the roles of $\psi$\ and $\phi$. 

Then by adding the same constant to $\phi$ and $\psi$ which 
obviously affects nothing, we can assume $0\leqslant\phi
\leqslant a$ where "$a$'' is a positive constant only 
depending on $A$ from the boundedness result before. 

Of course we can take a larger "$a$'', which we shall 
actually do below, as long as the dependence on $A$ is 
clear, or say finally we can still fix it to be some 
positive constant only dependent on $A$.   

As $lim_{_{t\to 0}}\gamma(t)=0$ by definition and the property 
of the function $\kappa$, we can fix $0<t_{0}<1$ sufficiently 
small such that $\gamma(t_{0}){t_{0}}^{n+3}<\frac{1}{3}$, which 
will also hold for $0<t<t_{0}$ since $\gamma$ is obviously 
decreasing. 

Fix such a $t$ for now and set $E_{k}=\{\psi<\phi-k a t\}$ 
where the "$a$'' is from above, but we still have not 
made the choice yet.

Clearly we have: 
$$\int_{E_{0}}g{\omega}^{n}=\frac{1}{2}\int_{E_{0}}\bigl((f+g)+
(g-f)\bigr){\omega}^{n}\leqslant\frac{1}{2}(1+\frac{1}{3})=\frac
{2}{3}.$$

Now we construct a function $g_{1}$ which is equal to $\frac{3g}
{2}$ over $E_{0}$ and some other nonnegative constant  for the complement. By the above estimate, it is easy to see 
that one can choose a proper constant (in [0,1]) such that $g_{1}
$ is still non-negative with $L^{p}$-norm bounded by $\frac{3A}{2}
$, and more importantly it has the proper total integral over $X$. 

So we can find a continuous solution $\rho\in PSH_{\omega}(X)$ 
as before by the approximation method such that 
$${\omega_{\rho}}^{n}=g_{1}{\omega}^{n}, ~~~max_{X}\rho=0$$ 
with lower bound of $\rho$ only dependent on $A$.  
\footnote{Notice we've used the existence of 
continuous solution at this point for the solution 
$\rho$.} By enlarging "$a$'' if necessary which clearly 
won't affect the set $E_{0}$, we can assume the 
lower bound of $\rho$ is $-a$. Now we can finally fix 
our constant "$a$'', and it clearly depends only on 
$A$ in an explicit way.

By noticing that $-2at\leqslant -t\phi+t\rho\leqslant 0$, 
it is easy to see 
$$E_{2}\subset E:=\{\psi<(1-t)\phi+t\rho\}\subset E_{0}.$$

Let's denote the set $\{f<(1-t^{2})g\}$ by $G$. Then 
over $E_{0}\setminus G$, we have: 
$$\bigl((1-t^{2})^{-\frac{1}{n}}\omega_{\phi}\bigr)^{n}\geqslant 
g{\omega}^{n}, ~~~\bigl((\frac{3}{2})^{-\frac{1}{n}}\omega_{\rho}
\bigr)^{n}=g{\omega}^{n}.$$

Hence we can conclude, using an inequality for mixed Monge-Amp\`ere measures from \cite{koj}, that over $E_0\setminus G$, 
$$(\frac{3}{2})^{-\frac{n-k}{n}}(1-t^{2})^{-\frac{k}{n}}{\omega_
{\phi}}^{k}\wedge{\omega_{\rho}}^{n-k}\geqslant g{\omega}^{n}.$$

\begin{rema}

This is a rather trivial result in smooth case which is just 
a direct application of arithmetic-geometric mean value 
inequality. Then by approximation argument, it should 
also hold in our case here. For the conclusion above, 
there is no need to restrict ourselves to the set $E_{0}\setminus G$. We 
can work globally on $X$ and use $g\chi_{_{E_{0}\setminus G}}{\omega}^n$ 
for the right hand side. 

Actually the rigorous approximation argument is local 
and uses nontrivial results about Dirichlet problem 
for Monge-Amp\`ere equation. The continuity of the 
functions is very involved in the proof which seems to 
be the main obstacle to carry over the whole argument 
in this part for merely bounded solutions. 

This is the point where the recent result in \cite{dinew} 
is applied. 

\end{rema}

Let's set $q=(\frac{3}{2})^{\frac{1}{n}}>1$, and rewrite 
the above inequality as:
$${\omega_{\phi}^{k}\wedge{\omega_{\rho}}^{n-k}}\geqslant q^{n-k}
(1-t^{2})^{\frac{k}{n}} g{\omega}^{n}$$ 
over $E_{0}\setminus G$. Now the following computation is quite 
obvious: \footnote{$t$ below can be taken to be sufficiently 
small, say $t<\frac{q-1}2$.}
\begin{equation}
\begin{split}
{\omega_{t\rho+(1-t)\phi}}^{n} 
&\geqslant \bigl((1-t)(1-t^{2})^{\frac{1}{n}}+q t\bigr)^{n}g{\omega}
^{n}\\
&\geqslant \bigl((1-t)(1-t^{2})+qt\bigr)^{n}g{\omega}^{n}\\
&\geqslant \bigl(1+t(q-1)-t^{2}\bigr)g{\omega}^{n}\\
&\geqslant \bigl(1+\frac{t}{2}(q-1)\bigr)g{\omega}^{n}.
\end{split}
\end{equation}

From the definition of $G$ and assumption of the theorem, 
we also have: 
$$t^{2}\int_{G}g{\omega}^{n}\leqslant\int_{G}(g-f){\omega}^{n}\leqslant
\gamma (t)t^{n+3}$$
which is just: 
\begin{equation}
\int_{G}g{\omega}^{n}\leqslant\gamma (t)t^{n+1}.    
\end{equation}
   
Hence we can have the following inequalities:  
\begin{equation}
\begin{split}
\bigl(1+\frac{t}{2}(q-1)\bigr)\int_{E\setminus G}g{\omega}^{n} 
&\leqslant \int_{E}{\omega_{t\rho+(1-t)\phi}}^{n}~~~~~~~~~~(the
~measure~inequality~(2.1))\\
&\leqslant \int_{E}{\omega_{\psi}}^{n}~~~~~~~~~~~~~~~~~~~(
comparison~principle)\\
&\leqslant \int_{E\setminus G}g{\omega}^{n}+\gamma (t)t^{n+1}~~
(the~integration~inequality~(2.2)).\nonumber
\end{split}
\end{equation}
and arrive at: 
$$\frac{q-1}{2}\int_{E\setminus G}g{\omega}^{n}\leqslant\gamma 
(t)t^{n}.$$

Therefore by noticing $E_{2}\subset E$, we get:
$$\frac{q-1}{2}(\int_{E_{2}}g{\omega}^{n}-\gamma (t)t^{n+1})
\leqslant\frac{q-1}{2}(\int_{E_{2}}g{\omega}^{n}-\int_{G}g{
\omega}^{n})\leqslant\frac{q-1}{2}\int_{E\setminus G}g{\omega
}^{n}\leqslant\gamma(t)t^{n},$$
and so we have 
$$\int_{E_{2}}g{\omega}^{n}\leqslant (t+\frac{2}{q-1})\gamma 
(t)t^{n}\leqslant\frac{3}{q-1}\gamma (t)t^{n}$$ 
for $t$ small enough. 

The claim proved before tells us: 
$$Cap_{\omega}(E_{4})\leqslant (\frac{a+1}{2at})^{n}\int_{E_{2}}
g{\omega}^{n}.$$

Combining this with the previous inequality, we have:
$$Cap_{\omega}(E_{4})\leqslant (\frac{a+1}{2a})^{n}\frac{3}{q-1}
\gamma (t).$$

Thus if $E':=\{\psi<\phi-(4a+2)t\}$ is nonempty, by the argument 
for boundedness result before, we should have:
$$2t\leqslant\kappa (Cap_{\omega}(E_{4}))\leqslant\kappa((\frac
{a+1}{2a})^{n}\frac{3}{q-1}\gamma (t))=t.$$
Clearly this is a contradiction for $t>0$. 

Anyway, we have from above that $\psi\geqslant\phi-(4a+2)t$. 

Hence $max_{X}(\psi-\phi)=max_{X}(\phi-\psi)\leqslant (4a+2)t$, 
which will give the desired conclusion.

\end{proof}

Now from this stability result, it is easy to get uniqueness 
result for continuous plurisubharmonic solutions after 
normalization.  

One can easily see the proof can be simplified a little if we 
only care about the uniqueness result. But this result above 
actually gives much better description of the variation of the 
solution under the perturbation of the right hand side of the 
equation (i.e, the measure).  

Now in the same vein as in \cite{koj} one gets the following corollary:

\begin{coro}\label{exp}
For any $\epsilon>0$, there exists $c=c(\epsilon,p,c_0)$,\ ($c_0$\ is an upper bound for $L^p$\ norms of $f$\ and $g$) such that
$$||\phi-\psi||_{\infty}\leq c||f-g||_1^{\frac{1}{n+3+\epsilon}}$$
provided $\phi$\ and $\psi$\ are normalised as before.
\end{coro}

Before we proceed further we make a small improvement of the stability exponent in the last corollary.
 
 Note that in the definition of set $G=\{f<(1-t^{2})g\}$\ one can exchange $t^2$\ with $\frac{t}{b}$\ for a sufficiently big independent constant $b$, and the the same argument still goes through, so $||f-g||_1\leq \gamma(t)t^{n+2}$\ implies $||\phi-\psi||_{\infty}\leq Ct.$ In particular the result in Corollary \ref{exp} holds with exponent $\frac1{n+2+\epsilon}$.

\section{Adjustment to Our Degenerate Case} 

Now we begin to adjust Ko\l odziej's argument for the 
situation in our main theorem. All the places which 
need to be considered have been pointed out at the 
spot. Let us now treat them one by one.   

\subsection{Comparison Principle}

In \cite{blo-koj}, authors constructed decreasing 
smooth approximation for bounded functions 
plurisubharmonic with respect to a K\"ahler 
metric. Using this, they got the following version 
of Comparison Principle,

\begin{theo}

For $\phi, \psi\in PSH_{\omega }(X)\cap L^\infty (X)$, where $(X, \omega)$ is 
a closed K\"ahler manifold, one has
$$\int_{\{\phi<\psi\} }(\omega+\sqrt{-1}\partial\bar{\partial }\psi)^n\leqslant\int_
{\{\phi<\psi\} }(\omega+\sqrt{-1}\partial\bar{\partial }\phi)^n.$$

\end{theo}

Though the result we want would be for 
some backround form $\omega\geqslant 0$, 
it would follow from the version above 
as we can perturb it by $\epsilon\omega_0$ with $
\omega_0>0$ and the constant $\epsilon>0$, since 
$X$ is K\"ahler. Those functions 
plurisubharmonic with respect to $\omega$ 
would still be plurisubharmonic with 
respect to $\omega+\epsilon\omega_0$. Using the
comparison principle above and letting $\epsilon\to 0$, 
we get the following version

\begin{theo}

For $\phi, \psi\in PSH_{\omega }(X)\cap L^\infty (X)$, where 
$X$ is a closed K\"ahler manifold and 
$\omega\geqslant 0$ is a real smooth $(1,1)
$-form over $X$, one has 
$$\int_{\{\phi<\psi\} }(\omega+\sqrt{-1}\partial\bar{\partial }\psi)^n
\leqslant\int_{\{\phi<\psi\} }(\omega+\sqrt{-1}\partial\bar{\partial }\phi
)^n.$$

\end{theo}
  
This is Comparison Principle for the 
adjusted argument for stability. 
 
\subsection{Inequalities for Mixed Measures} 
Our first observation is that although we considered our equations of the form
 $$\omega_{\psi}^n= f\omega^n,\ \ \omega_{\phi}^n= g\omega^n,$$
the volume form $\om^n$ played no significant role in the proof.
The only delicate point is the following inequality:

Suppose $\phi$\ and $\psi$ are continuous $\omega$-$PSH$ functions, and $f,\ g$\ are integrable functions on $X$. Suppose we have (locally or globally) the inequalities
$$ \omega_{\psi}^n\geq f\omega^n,\ \ \omega_{\phi}^n\geq g\omega^n,$$
then (locally where we have those inequalities or globally)
$$\forall_{ k\in\lbrace0, 1,\cdots, n\rbrace}\ \  \omega_{\psi}^k\wedge \omega_{\phi}^{n-k}\geq f^{\frac {k}{n}}g^{\frac{n-k}{n}}\om^n.$$
In other words we want to generalize the above inequality, for more general measures and moreover for {\it bounded} (i.e. not necessarily continuous) functions $\phi$\ and $\psi$. The following theorem is essentially taken from \cite{dinew}:
\begin{theo} Suppose the nonnegative Borel measure $\Omega$\ is well dominated by capacity, and let $\phi$\ and $\psi$\ be two bounded $\om$-psh functions on a K\"ahler manifold. Suppose the following inequalities hold
$$ \omega_{\psi}^n\geq f\Omega,\ \ \omega_{\phi}^n\geq g\Omega,$$
for some $f,\ g\in L^p(\Om),\ p>1$. Then
$$\forall_{ k\in\lbrace0, 1,\cdots, n\rbrace}\ \  \omega_{\psi}^k\wedge \omega_{\phi}^{n-k}\geq f^{\frac {k}{n}}g^{\frac{n-k}{n}}\Omega.$$
\end{theo}

In \cite{koj} (Lemma 1.2) this inequality was proved under the assumption that both $\phi$\ and $\psi$ are continuous and $\Omega=\om^n$. The proof is local, it can be rephrased in a setting in a ball in $\CC^n$. Then the argument goes via approximation for which a solution for the Dirichlet problem with boundary data is used. Since we deal with merely bounded functions (uppersemicontinuous by the plurisubharmonicity assumption), one cannot expect continuity on the boundary of the ball in general. But as observed in \cite{dinew} we can line-by-line follow the approximation arguments from \cite{koj} whenever the measure on the right hand side is the Lebesgue measure. Indeed, approximants at the boundary will not converge uniformly towards discontinuous boundary data, but the sequence of approximate solutions is again decreasing. This implies convergence in capacity by \cite{bed-tay}, which is enough for the argument to go through. In the case when $\om^n$\ is exchanged with a general measure well dominated in capacity one cannot rely only on the argument from \cite{koj}. But domination by capacity forces the measure $\Om$\  to vanish on pluripolar sets, hence one can use the result form \cite{dinew} to conclude. We refer to \cite{dinew} for the details. 

\section{Improvement on the Stability Exponent}

The exponent from Corollary \ref{exp} is quite important. In particular, since this inequality can be used to prove H\"older continuity for solutions of Monge-Amp\`ere equations with right hand side in $L^p$\ (see \cite{koj06}), the bigger the exponent in the inequality, the better H\"older exponent one can get. 

Trying to improve the exponent, one has to follow the main steps of the original proof and improve points where there is an exponent loss. Our strategy will be to iterate the original argument, defining at each step new function $\rho$\ and use the previous step to get estimates for $||\rho-\psi||_{\infty}$, which in turn will be used to choose the new set $E$\ in a ''better'' way.\\

The argument is divided into the following three parts.
 
The first part is the original argument quoted before with the improvement mentioned after Corollary \ref{exp}, 
which is the starting point for us. In the sequel the original argument will be often denoted as Step 1.

The second part, (i.e. Step 2), is the description of the 
iteration procedure. Since Step 1 differs slightly from all the others, we outline Step 2 below and sketch how to proceed throughout the next iterations. 

The mechanism is based on the fact that 
$||f-g||_1\leq\gamma(t)t^{\beta}$\ (in the improved original proof $\beta=n+2$) 
yields $\int_{\lbr\psi+kt<\phi\rbr}(\om+\ddc\psi)^n\leq c_0t^n$
for some constant $k$\ and $c_0$\ (in what follows $c_i$\ denote constants independent of the relevant quantities). So we try to find $\beta$\ as small as possible for which this implication holds true with uniform control on $c_0$ and enlarging $k$\ if needed. Note that from now on instead of $\om^n$\ we use the measure $\Omega$. It follows from the discussion above that Step $1$\ is not affected by that.

So assume $||f-g||_1\leq\gamma(t)t^{\beta},\ t<1$. Then if $l:=t^{\frac{\beta}{n+2}},\ \beta<n+2$, we obtain $||f-g||_1\leq\gamma(t)l^{n+2}$, so from Step $1$ we know that 
\begin{equation}\label{42}
\int_{E_2}g\Omega\leq\gamma(t)l^n,
\end{equation}
Where, as before $E_k:=\{\psi<\phi-kat\}$.
(Indeed, in Step $1$\ we have $t=l$, but one can check that the proof can be repeated in this situation). Hence
\begin{equation}\label{43}
\int_{E_2}g\Omega\leq c_1t^{\frac{\beta n}{n+2}},\ t\leq t_0
\end{equation}
(recall $\gamma(t)$\ decreases to $0$, as $t\searrow 0$).

Now fix a small positive constant $\delta$\ to be choosen later on.

Consider the ''new'' function 
$$g_1(z)=\begin{cases} (1+\frac{t^{\delta}}{2})g(z),\ &z\in E_2\\c_2g(z),\ &z\in X\setminus E_2,\end{cases}$$
where $0\leq c_2\leq 1$\ is choosen such that $\int_Xg_1\Omega=1$. (The constant $\frac 12$\ is taken to assure that the integral over $E_2$\ is less than $1$. Note that despite the fact that the  case $t$\ being small is of main interest, when $\delta$\ is also small the quantity $t^{\delta}$\ cannot be controlled by a constant smaller then $1$). As in Step $1$ we find a solution $\rho$ to the problem $(\om_{\rho})^n=g_1\om^n,\ max_X\rho=0.$ Again $\rho\geq-a$\ and we renormalize $\rho$\ by adding a constant so that $max_X(\psi-\rho)=max_X(\rho-\psi)$\ (this can by done in an uniform way).

Now by Step 1
\begin{align*}
&||\rho-\psi||_{\infty}\leq c_3||g-g_1||_1^{\frac 1{n+2+\epsilon}}=c_3([\int_{E_2}+\int_{X\setminus E_2}]|g-g_1|\Om)^{\frac 1{n+2+\epsilon}}=\\
=&c_3(2t^{\delta}\int_{E_2}g\Om)^{\frac 1{n+2+\epsilon}}\leq c_4t^{\frac{\delta+\frac{\beta n}{n+2}}{n+2+\epsilon}}.
\end{align*}
If $\delta$\ is sufficiently small the last exponent is less than $1$\ and we define $\al:=1-\frac{\delta+\frac{\beta n}{n+2}}{n+2+\epsilon}$.
Then by the above estimate
\begin{equation}
E_s=\lbr\psi+sat<\phi\rbr=\lbr(1-t^{\al})(\psi+sat)<(1-t^{\al})\phi
\rbr\subset
\end{equation}
\begin{align*}
&\subset\lbr\psi
<(1-t^{\al})\phi+t^{\al}\rho+c_4t-sat(1-t^{\al})\rbr=E\subset\\
&\subset\lbr\psi
<(1-t^{\al})\phi+t^{\al}\psi+2c_4t-sat(1-t^{\al})\rbr=\\
&=\lbr\psi+(sa-\frac{2c_4}{1-t^{\al}})t<\phi\rbr\subset E_k,
\end{align*}
provided $s\geq 4c_4+k$, (we take  $t<\frac 12$).

Consider the ''new'' set 
$$G_1:=\lbr f<(1-\frac{t^{\al+3\delta}}{8n2^{\frac{n-1}{n}}})g\rbr.$$
 Using  that $h(t)=(1+\frac{t^{\delta}}2)^{\frac 1n}-1-\frac1{4n2^{\frac{n-1}{n}}}t^{2\delta}$\ is increasing in $[0,1]$\ and hence nonnegative there, we conclude as in Step 1 that on $E_k\setminus G$
\begin{equation}\label{44}
(\om_{t^{\al}\rho+(1-t^{\al})\phi})^n\geq((1-t^{\al})
(1-\frac{t^{\al+3\delta}}{8n2^{\frac{n-1}{n}}})^{\frac 1n}+(1+\frac{t^{\delta}}{2})^{\frac 1 n}t^{\al})^ng\Omega\geq
\end{equation}
\begin{align*}
&\geq ((1-t^{\al})
(1-\frac{t^{\al+3\delta}}{8n2^{\frac{n-1}{n}}})
+(1+\frac1{4n2^{\frac{n-1}{n}}}t^{2\delta})t^{\al})^ng\Omega\geq
(1+\frac{t^{\al+2\delta}}{8n2^{\frac{n-1}{n}}})g\Omega.
\end{align*}

As in Step 1 on $G$\ we have
\begin{equation}\label{setG}
\frac{t^{\al+3\delta}}{8n2^{\frac{n-1}{n}}}\int_Gg\Omega\leq\int_G(g-f)
\Omega\leq\gamma(t)t^{\beta},
\end{equation}
so, using (\ref{44}), (\ref{setG}) and the comparison principle we obtain
\begin{equation}
(1+\frac{t^{\al+2\delta}}{8n2^{\frac{n-1}{n}}})\int_{E_k\setminus G}g\Omega\leq\int_E(\om_{(1-t^{\al})\phi+t^{\al}\rho})^n\leq
\end{equation}
\begin{equation*}
\leq\int_{E_k}g\Omega\leq\int_{E_k\setminus G}g\Omega+c_5\gamma(t)t^{\beta-\al-3\delta}.
\end{equation*}

Finally, as in Step 1, we obtain
$$\int_{E_k\setminus G}g\Omega\leq c_6\gamma(t)t^{\beta-2\al-5\delta}$$
and
$$\int_{E_s}g\Omega\leq c_7\gamma(t)t^{\beta-2\al-5\delta}.$$
If $\beta-2\al-5\delta=n$, we can proceed as in Step $1$\ to get $\max(\phi-\psi)=\max(\psi-\phi)\leq(2s+2)t$, and $||\phi-\psi||_{\infty}\leq C(\epsilon)||f-g||_1^{\frac{1}{\beta+\epsilon}},\ \forall\epsilon>0.$
Now $\beta-2\al-5\delta=n$\ yields
$$\beta(1+\frac{\frac{2n}{n+2}}{n+2+\epsilon})=n+2+5
\delta-\frac{2\delta}{n+2+\epsilon}.$$
It is clear that if $\delta$\ is sufficiently small $\beta$\ is smaller than $n+2$, hence we get an improvement.

Now in the last part we iterate the argument.

Consider $||f-g||_1\leq\gamma(t)t^{\beta_{k+1}}$, then as before $l=t^{\frac{\beta_{k+1}}{\beta_k}}$, $\int_{E_r}g\Omega\leq Ct^{\frac{n\beta_{k+1}}{\beta_k}}$, (compare with (\ref{42}), $r$\ is now chosen so that we can use the estimate on appropriate sublevel set from the previous step).

Choosing $\delta_{k+1}$\ small enough and proceeding in the same way as in the previous step one gets
$$\beta_{k+1}=n+2\al_{k+1}+5\delta_{k+1}.$$
($\al_{k+1}=1-\frac{\delta_{k+1}+\frac{\beta_{k+1} n}{n+2}}{n+2+\epsilon}$). This yields
\begin{equation}\label{requr}
\beta_{k+1}(1+\frac{2n}{\beta_{k}(\beta_{k}+\epsilon)})=
n+2+5\delta_{k+1}-2\frac{\delta_{k+1}}{\beta_{k}+\epsilon}
\end{equation}
If we choose $\lbrace\delta_{k}\rbrace$\ to be a sequence of sufficiently small numbers decreasing to $0$, one can obtain that $\lbrace\beta_{k}\rbrace$\ is decreasing (recall $n\geq 2$). If $A$\ is the limit of the sequence $\lbrace\beta_{k}\rbrace$\ one gets
$$A(1+\frac{2n}{A(A+\epsilon)})=n+2\Rightarrow A=\frac{n+2-\epsilon+\sqrt{(n-2-\epsilon)^2+8\epsilon}}{2}$$
Now $\epsilon\rightarrow0^{+}\Rightarrow A\rightarrow n$, so $\beta_k$'s can be arbitrarily close to n for $k$\ big enough if we take small enough $\epsilon$.

Thus this argument yields in paritcular Corollary \ref{exp2}.
\begin{rema}\label{rem4.1}
In the case when the measure $\Omega$\ is well dominated by capacity for $L^p$\ functions but the constant $\chi$\ is fixed one can construct $Q(t)$\ and aftrewards $\kappa(t),\ \gamma(t)$\ in such a way that $\gamma(t)\approx t^{\frac n{\chi}}$. Then one can use the same iteration technique as obove with  the exception that inequality (\ref{43}) should be improved to
$$\int_{E_2}g\Om\leq Ct^{\frac n{\chi}+\frac{\beta n}{n+2}}$$
(the factor $t^{\frac n{\chi}}$\ comes from the estimate of $\gamma$). The recurrence (\ref{requr}) now reads
\begin{equation}
\beta_{k+1}(1+\frac{2n}{\beta_{k}(\beta_{k}+\frac n{\chi})})=
n+2-\frac{\frac{n}{\chi}}{n+\frac{n}{\chi}}+5\delta_{k+1}-2\frac{\delta_{k+1}}{\beta_{k}+\frac n{\chi}}
\end{equation}
Again this is a convergent sequence and it can be computed that
$$\lim_{k\rightarrow\infty}\beta_k=n.$$
Hence the stability estimate in this case reads
\begin{equation}
||\phi-\psi||_{\infty}\leq c(\epsilon,c_0,X,\mu)||f-g||_{L^1(d\mu)}^{\frac 1{n+\frac n{\chi}+\epsilon}}
\end{equation}

\end{rema}
The following example shows that the exponent we obtained is sharp:
\begin{exam}\label{exa}
Fix appropriate positive constants $B,\ D$\ such that $D<B$\ and $B2^{2\al}<\log2+ D$, for some fixed $\al\in(0,1)$\ (such constants clearly exist). Then the function
$$\widehat{\rho}(z):=\begin{cases} B||z||^{2\al},\ &||z||\leq 1\\
\max\{ B||z||^{2\al},\log(||z||)+D\},\ &1\leq||z||\leq 2\\
\log(||z||)+D,\ &||z||\geq2\end{cases}
$$
is well defined, plurisubharmonic in $\CC^n$\ and of logarithmic growth. One can smooth out $\widehat{\rho}$, so that the new function $\rho$\ is again of logarithmic growth, radial, smooth away from the origin and $\rho(z)=B||z||^{2\al}$\ for $||z||\leq\frac 34$.

Via the standard inclusion
 $$\CC^n\ni z\longrightarrow[1:z]\in\mathbb P^n $$
one identifies $\rho(z)$\ with $$\overline{\rho}([z_0:z_1:,\cdots,:z_n]):
=\rho(\frac{z_1}{z_0},\cdots,\frac{z_n}{z_0})-\frac12
\log(1+\frac{||z||^2}{|z_0|^2})\in PSH(\mathbb P^n,\om_{FS})$$
(here $\om_{FS}$\ is the Fubini-Study metric on $\mathbb P^n$, and the values of $\overline{\rho}$\ on the hypersurface $\{z_0=0\}$\ are understood as limits of values of $\overline{\rho}$\ when $z_0$\ approaches $0$.) It is clear that $\om_{\overline{\rho}}^n=(dd^c\rho)^n$\ in the chart $z_0\neq 0$\ and in fact one can neglect what happens on the hypersurface at infinity. 

Now for a vector $h\in\CC^n$\ one can define $\rho_h(z):=\rho(z+h)$\ and analogously the corresponing $\overline{\rho}_h$. Note that when $||h||\rightarrow0,\ \overline{\rho}_h\rightrightarrows\overline{\rho}$.

One sees that
\begin{equation}
\label{al} B||h||^{2\al}\leq||\overline{\rho}_h-\overline{\rho}||_{\infty}
\end{equation}

The Monge-Amp\`ere measures of $\overline{\rho}$\ and $\overline{\rho}_h$\ are smooth functions except at the origin, and belong to $L^p(\om_{FS}^n)$, for some $p>1$\ dependent on $\al$.

Now $\int_{\mathbb P^n}|\om_{\overline{\rho}}^n-\om_{\overline{\rho}_h}^n|
=\int_{\CC^n}|(dd^c\rho)^n-(dd^c\rho_h)^n|$
To estimate the last term we divide $\CC^n$\ into three pieces (we suppose $||h||$\ is small):
$$\int_{\CC^n}|(dd^c\rho)^n-(dd^c\rho_h)^n|=
\int_{\{||z||\leq2||h||\}}+\int_{\{2||h||<||z||\leq \frac12\}}+\int_{\{||z||>\frac12\}}$$
Using the fact that $\rho$\ and $\rho_h$\ are smooth functions in a neighbourhood of $\{||z||>\frac12\}$\ one can easily estimate the last term by $||h||C_0$\ for some constant independent of $h$. For the first two terms we observe that $(dd^c\rho)^n=B^n||z||^{2n(\al-1)},\ (dd^c\rho_h)^n=B^n||z+h||^{2n(\al-1)}$.

Now we use a computation trick we found in \cite{KW}.
\begin{align*}
&\int_{\{||z||\leq2||h||\}}|(dd^c\rho)^n-(dd^c\rho_h)^n|
=\\
&=B^n\int_{\{||z||\leq2||h||\}}|||z||^{2n(\al-1)}-
||z+h||^{2n(\al-1)}|\leq\\
&\leq 2B^n\int_{\{||z||\leq3||h||\}}||z||^{2n(\al-1)}=C_1||h||^{2n\al}
\end{align*}
For the second term
\begin{align*}
&\int_{\{2||h||\leq||z||\leq\frac12\}}|(dd^c\rho)^n-(dd^c\rho_h)^n|=\\
&=B^n\int_{\{2||h||\leq||z||\leq\frac12\}}|||z||^{2n(\al-1)}-
||z+h||^{2n(\al-1)}|\leq \\
&\leq B^n\int_{2||h||<||z||}|\int_0^1<\nabla||z+th||^{2n(\al-1)},h>dt|\leq\\
&\leq C_2||h||\int_{||h||<||z||}||z||^{2n(\al-1)-1}\leq C_3||h||^{2n\al}
\end{align*}
provided $\al<\frac 1{2n}$, so that the integral is finite.
Finally we obtain for small $||h||$
\begin{equation}\label{cont}
\int_{\mathbb P^n}|\om_{\overline{\rho}}^n-\om_{\overline{\rho}_h}^n|\leq C_1||h||^{2n\al}+C_3||h||\leq C_4||h||^{2n\al}
\end{equation}
Suppose finally that we have a stability estimate $||\phi-\psi||_{\infty}\leq C_5||f-g||_1^{\frac 1 m}$. Then coupling
\ref{al} and \ref{cont} one gets
$$||h||^{2\al}\leq C_6(||h||^{2n\al})^{\frac 1 m},\ \al\in(0,\frac 1{2n})$$
If we let $||h||\rightarrow0$\ this can hold only if $m\geq n$.
\end{exam}
\begin{rema}
In \cite{EGZ} Authors show a stability estimate of another type: In the setting as above ($\Om$\ is now equal to $\om^n$) 
\begin{equation}
||\phi-\psi||_{\infty}\leq c(\epsilon,c_0,\omega)
||\phi-\psi||_{L^2(\om^n)}^{\frac2{nq+2+\epsilon}}
\end{equation}
($c_0$\ is a constant that controls $L^p$\ norms of Monge-Amp\`ere measures of $\phi$\ and $\psi$). Using the same reasoning as in \cite{EGZ}
one can show more generally that
\begin{equation}
||\phi-\psi||_{\infty}\leq c(\epsilon,c_0,\omega)
||\phi-\psi||_{L^s(\om^n)}^{\frac s{nq+s+\epsilon}}, \forall s>0.
\end{equation}
Using the same example and similar estimates one can show that this exponent is also sharp, provided that $p<2$\ and $s>\frac{2np}{2-p}$\ (the reason for these obstructions is that the second integral we estimate as in the example would be divergent otherwise). It is, however, very likely that these exponents are sharp in general.
\end{rema}

\section{Continuity of Solutions in the Case of a Pullback Form via a Locally Birational Map}
\label{continuity}

We give below a more detailed proof of the continuity statememt in Theorem \ref{1.1}. 
Arguments used heavily rely on \cite{koj98} and at some places we just follow it line 
by line. This section is unrelated with the other ones in the note. Recall once again, 
that this result is known already.

First of all we recall the geometrical background. Let $X$\ be the base closed K\"ahler manifold we work on, and $F:X\rightarrow \mathbb{CP}^N$, is a map with the property that the image $F(X)$\ has the same dimension and $F$\ is itself locally birational i.e. for every small enough neighbourhood $U$\ of any point on $F(X)$, each component of $F^{-1}(U)$\ is birational to $U$. A typical global example of this situation is obtained as follows: if $X$\ carries a big line bundle $L$, the linear series corresponding to $L^n$\ generate (for sufficiently big $n\in\mathbb N$) a birational morphism into $\mathbb{CP}^N$\ with the claimed properties. Note hovewer that local and global birationality are different notions (see the example below) and if one has to deal with the global birationality one has to impose some additional assumptions for the argument to go through.

Consider now $Y:=F(X)$. By the Proper Mapping Theorem $Y$\ is a (singular in general) subvariety in $\mathbb{CP}^N$. It is also clear that $Y$\ is irreducible and locally irreducible variety (the latter follows from the local birationality). Recall that an upper semicontinuous function $u$\ on a singular variety $W$\ is called weakly plurisubharmonic if for every holomorphic disc $f:\Delta\rightarrow W$\ the function  $u\circ f$ is a subharmonic function (see \cite{for-nar}). In that paper it is proved (in fact in a much more general situation of Stein spaces) that any such function $u$\ can be extended locally to the ambient space to a classical plurisubharmonic function i.e. for every $x\in Y$ there exists a small Euclidean ball $B$ in $\mathbb{CP}^N$, centered at $x$\ and a function $v\in PSH(B)$, such that $v|_{B\cap Y}=u$.

Now suppose $\phi$\ is a positive discontinuous solution of the Monge-Ampere equation in question and let $d:=\sup(\phi-\phi_{*})>0$, where $\phi_{*}$\ denotes the lower semicontinuous regularization of $\phi$. Note that the supremum is attained, and if $E$\ is the closed set $\lbrace\phi-\phi_{*}=d\rbrace$, there exists a point $x_0$\ such that $\phi(x_0)=\min_E\phi$. Positivity is a technical assumption that can always be achieved by adding appropriate constant since we already know that $\phi$\ is bounded.

By assumption there exist analytic sets $Z\subset X$\ and $W\subset Y=F(X)$\ such that $F|_{X\setminus Z}\rightarrow Y\setminus W$\ is a biholomorphism and moreover $S:=\lbrace \omega^n=0\rbrace\subset Z$. Note that in the general case of a big form $S$\ need not be contained in an analytic set- it may well happen that $S$\ is open in $X$. 

Two possibilities might take place
\begin{enumerate}
\item $x_0\in X\setminus S$. In this case $\omega$\ is strictly positive in a small ball centered at $x_0$\ and repeating the argument from Section 2.4 in \cite{koj98} we obtain a contradiction.
\item $x_0\in S$. Then we shall produce a domain $V$\ (not contained in a chart in general) and a potential $\theta$ of $\omega$\ in $V$\ with the property that $\inf_{\partial V}\theta>\theta(x_0)+b$, where $b$\ is a positive constant.
\end{enumerate}

Consider $F(x_0)=z$\ and a neighbourhood $U$\ of $z$\ in $Y$, such that its preimages are birational to  it. Choose the one $x_0$\ sits in. For the rest of the argument we restrict ourselves to $F|_{F^{-1}(U)\ni x_0}\rightarrow U$.  
Consider the pushforward function
$$F_{*}\phi(z):=\begin{cases}\phi(w),\ \ if\  z\in Y\setminus W, w\in X\setminus Z,\  F(w)=z\\
\limsup_{X\setminus Z\ni \zeta\rightarrow z}F_{*}(\zeta)
\end{cases}$$
and a local potential $\eta$\ for the K\"ahler form on $U\cap\mathbb {CP}^N$.

Claim: $\eta+F_{*}\phi$\ is weakly subharmonic on $Y$.

Proof: Weak subharmonicity is a local property, hence it is enough to check it in a neighbourhood 
of any point on $Y$. For regular points of $Y$\ this is evident. However at singularities of $Y$\ one might a priori run into trouble as the example of a double point shows. Indeed, consider the following (classical) local example:

Let $$F:\mathbb C\ni t\rightarrow (t^2-1,t(t^2-1))\in\mathbb C^2$$
The image $F(\mathbb C)$\ sits in the variety $\lbrace(z_1,z_2)\in\mathbb C^2|z_1^2+z_1^3=z_2^3\rbrace$. Observe that $F$\ is a bijection onto its image, except for the points $1$\ and $-1$\ being mapped to $(0,0)$. But then it is clear that the pushforward of a subharmonic function $w$\ on $\mathbb C$\ cannot be weakly subharmonic on the image if $w(1)\neq w(-1)$. Note that $F$\ is not locally birational though. 

Observe that local birationality forces the analytic set $Y$\ to be locally irreducible. 
Then there is a classical theorem (see \cite{De}, Theorem 1.7) stating that on a locally irreducible variety 
$Y$\ and a locally bounded plurisubharmonic function $w$\ defined on $\Reg Y$- the 
regular part of $Y$ the extension via limsup technique $w(z):=\limsup_{\zeta\rightarrow 
z,\ \zeta\in \Reg Y}w(\zeta)$\ is a weak plurisubharmonic function. Moreover, it follows from the proof that for any $s\in Y$\ and any birational modification $G:Y^{'}\rightarrow Y$\ of $Y$\ the pulled-back function $G^{*}w$\ is constant on the fiber $G^{-1}(s)$.  

Now if $\omega_M$\ is the K\"ahler metric which defines $\omega$\ (i.e. $\omega=F^{*}
\omega_M$), fix  $\rho$- a local potential of $\omega_M$\ 
near $z$\ (in $\mathbb{CP}^n$). First we shall modify $\rho$\ exactly as in \cite{koj98}:

In local coordinates in a ball $B^{''}$\ centered at $z$\ $\rho$\ is strictly plurisubharmonic 
smooth function and is expanded as 
\begin{equation}
\begin{split}
\rho(z+h)
&=  \rho(z)+2\Re (\sum_{j=1}^na_jh_j+\sum_{j,k=1}^nb_{jk}h_jh_k)+\sum_{j,k=1
}^nc_{j\bar{k}}h_j\bar{h}_k+o(|h|^2) \\
&= \Re P(h)+H(h)+o(|h|^2), \nonumber
\end{split}
\end{equation}
where $P$ is a complex polynomial in $h$\ and $H$\ is the compex Hessian at $z$.

Proceeding exactly as in \cite{koj98} (Lemma 2.3.1) $\eta:=\rho-\Re P(\cdot-z)$\ is also a local potential for $\omega_M$, with the additional property that $\eta$\ has a strict local miniumum at $z$\ (we use at this point that $H$\ is strictly positive definite). This means that for a smaller ball (which after possible shrinking we again denote by $B^{"}$) $\inf_{\partial B^{"}}\eta>\eta(z)+b^{"}$\ for some positive constant $b^{"}$. Adding a constant if necessary one can further assume that $\eta(z)>0$.

Now by Fornaess-Narasimhan theorem  we find a small euclidean ball $B^{'}$\ in $\mathbb{CP}^n$\ centered at $z$\ and a function $\psi\in PSH(B^{'})$, such that $\psi|_{Y\cap B^{'}}=\eta+\phi$. On a neighbourhood of a slightly smaller ball $B$\ (everything is contained in $B^{'}$\ and $B^{"}$) $\psi$\ can be approximateded by a sequence on smooth plurisubharmonic functions $\psi_j$\ decreasing towards it. Again (decreasing a bit $b^{"}$\ if necessary) one obtains  $\inf_{\partial B}\eta>\eta(z)+b$\ for some positive constant $b$.
 Now we pull back the ball and the regularizations: let $V:=F^{-1}(B^{"}\cap Y)$\ and $u_j:=\psi_j(F(w))$ ( $u_j$\ are assumed to be defined only on small neighbourhood of  $V$). Of course these are continuous plurisubharmonic functions on $V$\ which decrease towards $u:=\eta\circ F+F^{*}(F_{*}\phi)=\eta\circ F+\phi$ (the equality is due to the fact that $\phi$\ has to be a constant on the fiber). Note that $V$\ need not be an Euclidean domain anymore (i.e. it need not be contained in a coordinate chart), nevertheless $\eta\circ F$\ is a global potential of $\omega$\ on this set. This is the essential difference between this special situation and the general case.
 Next we state a lemma which is essentially contained in \cite{koj98} (Section 2.4). We include the proof for the sake of completeness.
\begin{lemm} There exist $a_0>0,\ t>1$ such that the sets
$$W(j,c):=\lbrace tu+d-a_0+c<u_j\rbrace$$
are non-empty and relatively compact in $V$\ for every constant $c$\ belonging to an interval which does not depend on $j>j_0$.
\end{lemm}
\begin{proof}
Note that $E(0):=\lbrace u-u_{*}=d\rbrace\cap \overline{V}=E\cap \overline{V}$, since the potential is continuous. Also the sets $E(a):=E:=\lbrace u-u_{*}\geq d-a\rbrace\cap \overline{V}$\ are closed and decrease towards $E(0)$. Hence if $c(a):=\phi(x_0)-\min_{E_a}\phi$\ we have that $\limsup_{a\rightarrow 0^{+}}c(a)\leq 0$, for othrwise we would get a contradiction with the definition of $d$. Hence 
$$c(a)<\frac13 b$$
for $0<a<a_0<min(\frac13 b, d)$. Let $A:=u(x_0)$. Note that $A>d$\ since the potential is greater than $0$\ at $x_0$, and $\phi$\ as a globally positive function has to be greater than $d$\ at $x_0$. One can choose $t>1$, such that it satisfies
$$(t-1)(A-d)<a_0<(t-1)(A-d+\frac 2 3 b)$$
Now if $y\in\partial V\cap E(a_0)$\ one gets
$$u_{*}(y)\geq \eta(F(x_0))+b+F^{*}F_{*}\phi(x_0)\geq A-d+\frac 23 b$$
Hence $u(y)\leq u_{*}(y)+d<tu_{*}(y)+d-a_0$. Note that this inequality extends to a neighbourhood  of $\partial V\cap E(a_0)$. Taking another neighbourhood relatively compact in the first and applying Hartogs type argument one obtains
$$u_j<tu(y)+d-a_0,\ \ \forall j>j_1$$
For the rest part of $\partial V$\ the same inequality holds if we take big enough $j_1$\
and the proof is even simpler, since $u-u_{*}$\ is less than $d-a_0$\ there. This proves the relative compactness on $W(j,c)$\ in $V$.

Note that from the left inequality defining $t$\ one gets $(t-1)u_{*}(x_0)<a_0$, hence $$tu_{*}(x_0)<u(x_0)-d-a_1+a_0<u_j(x_0)-d-a_1+a_0$$
for some constant $a_1>0$. This implies that the sets $W(j,c),\ c\in(0,a_1)$\ contain some points near $x_0$, hence they are non empty. 
\end{proof}

Now applying Lemma 2.3.1 from \cite{koj98} (one can verify that despite the fact that $V$\ can be a non Euclidean set the argument still goes through) one can bound the capacity $cap(W(j,a_1),V)$\ from below by an uniform positive constant. On the other hand $W(j,a_1)\subset\lbrace u+(d-a_0+a_1)<u_j\rbrace$ and this contradicts the fact that the decreasing sequence $u_j$\ has to converge towards u in capacity. This proves that $d=0$, hence $\phi$\ is continuous.
\begin{rema} As we have seen the argument cannot be applied in the case of a (globally) birational map. Then further assumption is needed to assure the pushforward to be plurisubharmonic. A satisfactory additional assumption is that the fibers in the preimage have to be connected. Then the function has to be constant on any nontrivial connected fiber and this is enough to push it forward onto the image.
\end{rema}
\section{Remarks} 

Complex Monge-Amp\`ere equations are of great 
interest in geometry. In \cite{thesis}, the following 
version of the Monge-Amp\`ere equation
$$(\omega+\sqrt{-1}\partial\bar{\partial}u)^n=e^u\Omega$$
is the main focus. Of course anything new would be for a degenerate class $[\omega]$\ as in the settings of Theorem 1.1. And using the argument in 
\cite{koj}, we know that the main result in this 
work would also apply for it. To be precise, the following folds:
\begin{theo} Let $\omega$\ be a big form and $u_1$,\ $u_2$\ be $\om$-psh solutions for the following Monge-Amp\`ere equations:
$$(\om_{u_1})^n=e^{u_1}\Om_1,\ \ (\om_{u_2})^n=e^{u_2}\Om_2,$$
where $\Om_1$\ and $\Om_2$\ are smooth volume forms such that
$$\int_X|\Om_1-\Om_2|\leq\gamma(t)t^{n+3}.$$
Then $$||u_1-u_2||_{\infty}\leq Ct.$$
\end{theo} 
\begin{proof}Since the comparison principle for big forms is avialiable the proof is entirely the same as in Theorem 5.2 in \cite{koj}.
\end{proof}
The following problems are related to the results in \cite{koj06} and \cite{EGZ}, stating that when $\om$\ is K\"ahler form on a compact K\"ahler manifold, the solutions of 
$$\om_{\phi}^n=f\om^n, f\in L^p(\om^n),\ p>1$$
 are H\"older continuous. In general the H\"older exponent depends on the manifold $X$, and on $n$\ and $p$ (\cite{koj06}). Under the additional assumption that X is {\it homogeneous} i.e. the automorphism group $\Aut(X)$\ acts transitively the exponent is independent of $X$ and is not less that $\frac 2{nq+2},\ q=\frac p{p-1}$\ (\cite{EGZ}). One can ask the following questions:
\begin{enumerate}
\item Is the solution continuous  when $\om$\ is semi-positive and big in general?
\item If so, does the H\"older continuity hold in the case $\om$\ is merely semi-positive and big? 
\item Does the H\"older exponent on general manifold do really depend on the manifold? In the corresponding result in the flat case (\cite{GKZ}) the H\"older exponent is uniform and independent of the domain. Moreover the proof in \cite{koj06} strongly depends on a regularization procedure for $\om$-psh functions, which consists of patching local regularizations, and this is the point where the geometry of the manifold influences the exponent. In particular are there another regularization procedures of more global nature that are not so affected by the local geometry?
\item Is the exponent for the homogeneous case sharp? Note that in the flat case in \cite{GKZ} there is also a gap between the exponent given there $\frac 2{qn+1}$\ and the exponent $\frac 2{qn}$, for which an example is shown.
\item It is interesting to compare the stability resuts we have proven and the one in \cite{EGZ}. In particular, is the stability exponent in \cite{EGZ} sharp in general?
\item It would be very interesting to generalize H\"older continuity to more singular measures. One possible application of such a result would be a criterion for H\"older continuity of the Siciak extremal function of certain compact sets in $\CC^n$\ (see \cite{kojnotes} for a definition). Such a property is very important from pluripotential point of view. So one has to study the equilibrium measure of the compact. The problem is that such measures are singular with respect to the Lebesgue measure, while \cite{koj06}\ and \cite{EGZ} rely strongly on smoothness of $\om^n$. However, as this note shows, some arguments can be adjusted to singular measures either.
\end{enumerate}
We hope to address some of these questions in the future.

\begin{flushleft}Jagiellonian University\\
Institute of Mathematics\\
Reymonta 4, 30-059 Krak\'ow, POLAND.\\
E-mail {\tt slawomir.dinew@im.uj.edu.pl}\\
\end{flushleft}

\begin{flushleft}Michigan University\\
Department of Mathematics\\
4835 East Hall, Ann Arbor, USA\\
E-mail {\tt zhangou@umich.edu}\\
\end{flushleft}

\bigskip
 Key words and phrases: K\"ahler manifold, complex Monge-Amp\`ere operator.\\
 2000 Mathematics Subject Classification: Primary 32U05, 53C55; secondary: 32U40. 


\begin{thebibliography}{$$}

\bibitem{BGZ} Benelkourchi, Slimane; Guedj, Vincent; Zeriahi, Ahmed: 
A priori estimates for weak solutions of complex Monge-Amp\`ere 
equations. Preprint. ArXiv 0704.0866.

\bibitem{bed-tay} Bedford, Eric; Taylor, B. A.: 
A new capacity for plurisubharmonic functions. 
Acta Math. 149 (1982), no. 1-2, 1--40.

\bibitem{blo-koj} B\l ocki, Zbigniew; Ko\l odziej, 
S\l awomir: On regularization of plurisubharmonic 
functions on manifolds. Proc. Amer. Math. Soc. 
135 (2007), no. 7, 2089--2093 (electronic). 

\bibitem{dinew} Dinew, S\l awomir: An Inequality 
for Mixed Monge-Amp\`ere Measures. Preprint. ArXiv 0705.0974.

\bibitem{De} Demailly, Jean-Pierre: Mesures de Monge-Amp\`ere et 
caracterisation g\'eom\'etrique des vari\'et\'es alg\'ebriques affines. M\'emoires 
Soc. Math. France 19 (1985) 1-125.

\bibitem{DP} Demailly, Jean-Pierre; Pali, Nefton: Degenerate complex 
Monge-Amp\`ere equations over compact K\"ahler manifolds. Preprint. 
ArXiv 0710.5109.

\bibitem{EGZ} Eyssidieux, Phillipe; Guedj, Vincent; 
Zeriahi, Ahmed: Singular K\"ahler-Einstein metrics. 
Preprint. ArXiv math/0603431.

\bibitem{for-nar} Fornaess, John Erik; 
Narasimhan, Raghavan: The Levi problem on 
complex spaces with singularities. 
Math. Ann. 248 (1980), no. 1, 47--72. 

\bibitem{GKZ} V. Guedj, S. Ko\l odziej, A. Zeriahi: H\"older continuous 
solutions to Monge-Amp\`ere equations. preprint. ArXiv math/0607314.

\bibitem{koj98} Ko\l odziej, S\l awomir: The complex 
Monge-Ampere equation. Acta Math. 180 (1998), 
no. 1, 69--117.

\bibitem{koj} Ko\l odziej, S\l awomir: The Monge-Ampere 
equation on Compact K\"ahler Manifold. Indiana Univ. 
Math. J. 52 (2003), 667-686.

\bibitem{kojnotes} Ko\l odziej, S\l awomir: The complex 
Monge-Ampere equation and pluripotential theory. 
Mem. Amer. Math. Soc. 178 (2005), no. 840, x+64 pp. 

\bibitem{koj06} Ko\l odziej, S\l awomir: H\"older continuity 
of solutions to the complex Monge-Ampere equation 
with the right hand side in $L^p$. Preprint. ArXiv math/0611051.

\bibitem{KW} H. Kozono, H. Wadade: Remarks on 
Gagliardo-Nirenberg type inequality with critical Sobolev 
space and BMO. to appear in Math Z.

\bibitem{t-znote} Tian, Gang; Zhang, Zhou: On the 
K\"ahler-Ricci flow on projective manifolds of general 
type. Chinese Annals of Mathematics - Series B, Volume 
27, Number 2, 179--192.   

\bibitem{zh} Zhang, Zhou: On Degenerate Monge-Ampere 
Equations over Closed K\"ahler Manifolds. Int. Math. Res. 
Not. 2006, Art. ID 63640, 18 pp. 

\bibitem{thesis} Zhang, Zhou: Degenerate Monge-Ampere 
Equations over Projective Manifolds. PHD Thesis at MIT, 
2006. 

\end{thebibliography}
\end{document}